\documentclass[12pt]{article}
\usepackage{srcltx}
\usepackage{amssymb,amsmath,amsfonts,amsthm}
\usepackage{cite}
\newtheorem{lem}{Lemma}

\begin{document}

UDC 512.542

\begin{center}\textbf{Recognizability by spectrum of alternating groups}\end{center}

\begin{center}Gorshkov I.B.\footnote{Supported by Russian Foundation for Basic Research (Grants~11-01-00456, 12-01-31221, 12-01-90006).}\end{center}

\hfill \textit{Dedicated to 70-th anniversary of V. D. Mazurov}

~

\textsl{Review}: The spectrum of a group is the set of its element orders. A finite group $G$ is said to be recognizable by spectrum if every finite group that has the same spectrum as $G$ is isomorphic to $G$. We prove that the simple alternating groups $A_n$
are recognizable by spectrum when $n\neq 6, 10$. This implies that every finite group with the same spectrum as that of a finite nonabelian simple group, has at most one nonabelian composition factor.

~


~

Suppose $G$ is a finite group, $\pi(G)$ is the set of prime divisors of its order, $\omega(G)$ is the spectrum of $G$, i.\,e. the set of its element orders. The Gruenberg-Kegel graph, or the prime graph, $GK(G)$ is defined as follows. The vertex set of the graph is $\pi(G)$. Two distinct primes $p$ and $q$ of $\pi(G)$ seen as verticies of the graph $GK(G)$, are connected by an edge if, and only if, $pq\in\omega(G)$. A group $G$ is said to be recognizable by spectrum (shortly, recognizable) if for every finite group $L$ the equality $\omega(L)=\omega(G)$ implies that $L\simeq G$.
In \cite{MazKon}, \cite{Zav} the following groups were proved to be recognizable: $A_{n}$, $n\geq5$, $n\not\in\{6,10\}$ and $n-p\in \{0,\ 1,\ 2\}$, where $p$ is a maximal prime number which does not exceed $n$. The proof is based on the fact that the graph $GK(A_{n})$ is disconnected, and the prime $p$ forms its connected component, which is not true in general case. The group $A_6$ was proved to be irrecognizable in \cite{Shi}. In \cite{Maz2} it was shown that $A_{10}$ is irrecognizable. The groups $A_{16}$ and $A_{22}$ were proved to be recognizable in \cite{Zav}, \cite{Changg} respectively. In particular, the question on recognizability is settled for all alternating groups $A_k$, where $5\leq k\leq 25$.

A finite simple group $L$ is said to be quasirecognizable, if every finite group $G$ with $\omega(L)=\omega(G)$  has the unique nonabelian composition factor, which is isomorphic to $L$. It was shown in \cite{ZavMaz} that quasirecognizability of a finite simple alternating group implies its recognizability. In \cite{Vac} a theorem describing the properties of chief series of groups with the same spectrum as that of the alternating group was proved. Using the above results we were able to prove that all simple alternating groups except for $A_6$ and $A_{10}$, are recognizable.

\textbf{Theorem 1. }\emph{
Suppose $G$ is a finite group with $\omega(G)=\omega(A_n)$, where $n\geq5,\ n\neq6,\ 10$. Then $G$ is isomorphic to $A_n$.
}

This theorem gives an affirmative answer to Question 16.107 in the Kourovka Notebook. Also together with available results on groups $A_6$ and $A_{10}$ it implies that the Question 16.27 is answered in the affirmative. Moreover, as it is indicated in the comment to Question 16.27, Theorem 1 and \cite[Corollary 7.3]{VasVd} entail the following statement.

\textbf{Theorem 2. }\emph{ Suppose $L$ is a finite nonabelian simple group and $G$ is a finite group with $\omega(G)=\omega(L)$. Then $G$ has at most one nonabelian composition factor. Moreover, if the group $L$ is distinct from the groups $\operatorname{PSL}_3(3)$, $\operatorname{PSU}_3(3)$, $\operatorname{PSp}_4(3)$, then $G$ has exactly one nonabelian composition factor.
}

Note that for every exception in Theorem 2 there exists a soluble group with the same spectrum.

\textbf{\S 1. Preliminaries}

\begin{lem}\label{factor}
Suppose $K$ is a normal subgroup of a finite group $G$, $\overline{G}=G/K$, $\overline{x}$ is the image of an element $x$ of $G$ in the group $\overline{G}$. If $(|x|,|K|)=1$, then $C_{\overline{G}}(\overline{x})=C_G(x)K/K$.
\end{lem}
\textsl{Proof.} See \cite[Theorem 1.6.2]{Khu}.\\

\begin{lem}\label{vas2} Assume that a finite group $G$ has a normal series $1 <K < M < G$, and pairwise distinct prime numbers $p,\ q$ and $r$ are such that $p$ divides $|K|$, $q$ divides
$|M/K|$ and $r$ divides $|G/M|$. Then the numbers $p$, $q$ and $r$ can not be pairwise non-adjacent in $GK(G)$.
\end{lem}
\textsl{Proof.} See \cite[Lemma 1.1]{Vas}.

\begin{lem}\label{vac2}
 Every odd number of $\pi(Out(P))$, where $P$ is a finite simple group, belongs to the spectrum of the group $P$ or does not exceed $m/2$, where $m=max_{p\in\pi(P)}p$.
\end{lem}
\textsl{Proof.} See \cite[Lemma 14]{Vac}.

\begin{lem}\label{zav}
If $n\geq5,\ n\neq6,\ 10$ and $n=r+m$, where $r$ is a prime which is greater than $3$, $m\in\{0,\ 1,\ 2\}$, then the group $A_n$ is recognizable.
\end{lem}
\textsl{Proof.} See \cite{MazKon} and \cite{Zav}.

\begin{lem}\label{AS}
If $n>6$, then $\omega(A_n)\not\subseteq \omega(S_{n-1})$.
\end{lem}
\textsl{Proof.} See\cite[Lemma 1.2]{Zav}
\begin{lem}\label{ASS}
If $n>6$, then $\omega(A_n)\neq \omega(S_{n})$.
\end{lem}
\textsl{Proof.} See\cite{Zav}

\begin{lem}\label{zav2}
Suppose that a quotient group $H = T/N$ of a finite group $T$
is isomorphic to a symmetric or alternating group of degree $m$, where $m \geq 5$ and $N\neq 1$. Then $\omega(T/N)\neq\omega(T)$.
\end{lem}

\textsl{Proof.} See \cite{ZavMaz}.

\begin{lem}\label{vac}
Suppose that $\omega(G)=\omega(A_n)$, $n>21$. Then a chief series of $G$ has a factor, which is isomorphic to $A_k$, for some $k$ of $[p,n]$, where $p$ is the maximal prime not exceeding $n$. In addition, no other factor contains $p$ in its spectrum.
\end{lem}

\textsl{Proof.} See \cite{Vac}.

\begin{lem}\label{sp}
Suppose that $g\in A_n$ is an element of order $r$, where $r$ is a prime and $n/2<r<n-2$, $n\geq 5$. Then $C_{A_n}(g)\simeq \langle g\rangle\times A_{n-r}$.
\end{lem}
\textsl{Proof.} An easy check.

\medskip
\textbf{\S 2. Proof of Theorem 1}
\medskip

Suppose that $n\geq 26$ is the least number such that $L= A_n$ is not recognizable, $G$ is a finite group with $\omega(G)=\omega(L)=\omega$ and $G\not\simeq L$, $1=G_0<G_1<...<G_t=G$ is a chief series of $G$, $R_i=G_i/G_{i-1}$, $i\in[1,t]$. According to Lemma \ref{vac} there exists $m$ such that $R=R_m\simeq A_k$, $k\in[p,n]$, where $p$ is the largest prime not exceeding $n$. We chose a chief series so that the quotient $G/G_{m-1}$ has the least order. Then the group $G/G_{m-1}$ possesses the unique minimal normal subgroup, which is isomorphic to $A_k$. So $G/G_{m-1}$ is isomorphic either to $A_k$, or to $S_k$. By Lemma \ref{zav} the number $n-p$ is greater than $2$. Put $\Pi=\{r|n/2<r<p\}$, $\overline{\Pi}=\Pi \setminus\pi(G_{m-1})$. Note that numbers of $\Pi$ are pairwise non-adjacent in $GK(L)$, and $\overline{\Pi}$ consists of those elements from $\Pi$, who divide just the order of the factor $R$. Let $\mu=\mu(G)$ denote the set of maximal under divisibility elements of $\omega(G)$. Note that $\omega(G)$
is uniquely determined by the set $\mu(G)$.

\begin{lem}\label{tc}
Under $n\neq 27$ there exist two distinct primes $r_1, r_2\in \Pi$ such that $n-r_i\not \in \{3,4,6,10\}$, $i=1,2$.
\end{lem}
\textsl{Proof.} Since $p>r$ for every $r\in\Pi$ and $n-p>2$, we have $n-r>4$. Having $n\geq42$, it is easy to check that $|\Pi|\geq4$ (see, for example, \cite[Lemma 1]{MazKon}). If $26\leq n<48, n\neq 27$, then the statement of the lemma can be checked directly. Under $n=27$ the statement is false.

\begin{lem}\label{centralizer}
Let $y\in G$ be an element of order $r$, where $r\in \Pi\cup\{p\}$. If a Sylow $r$-subgroup of $G$ is cyclic of prime order, then $\omega(C_G(y))=\omega(C_L(x))$, where $x$ is an $r$-element of the group $L$.
\end{lem}
 \textsl{Proof.} If an element $z\in G$ centralizes some $r$-element of $G$, then there exists its conjugate  $z'$, which centralizes $y$. Therefore, if $rm\in \omega$, then $rm\in \omega(C_G(y))$. The same we can say of the group $L$ and the element $x$. Thus the sets $\mu(C_G(y))$ and $\mu(C_{L}(x))$ consist precisely of those elements of the set $\mu$, who are divided by $r$.

\begin{lem}\label{centralizer2}
Suppose that $y\in G$ is an element of order $r$, where $r\in \Pi\cup\{p\}$, and a Sylow $r$-subgroup of $G$ is cyclic and has a prime order. Then $C_G(y)=\langle y\rangle\times M$, where $\mu(M)=\mu(Alt_{n-r})$.
\end{lem}
\textsl{Proof.} Put $M=C_G(y)/\langle y\rangle$. A Sylow $r$-subgroup of $G$ is cyclic and has prime order. Hence $\langle y\rangle.M\simeq \langle y\rangle\times M$ and $\mu(M)=\{l/r| l\in \mu(C_G(y))\}$. If $x$ is an $r$-element of $L$, then according to Lemma \ref{sp} we have $C_L(x)=\langle x\rangle\times N$, where $N\simeq Alt_{n-r}$. Therefore, $\mu(M)=\mu(N)=\mu(Alt_{n-r})$.

\begin{lem}\label{pi}
$|\overline{\Pi}|\geq |\Pi|-1$
\end{lem}
\textsl{Proof.}
 Assume that there exists $i$ such that $1\leq i< m$ and $|R_i|$ is divided by two numbers of $\Pi$. Suppose that $R_i=T_1\times T_2\times...\times T_h$, where $T_j\simeq T_1,\ 1<j\leq h$, $T_1$ is a finite simple group. The fact that for every $r_1,\ r_2\in \Pi,\ r_1\neq r_2$, the set $\omega$ does not contain a number $r_1r_2$ implies that $R_i=T_1$. By Lemma \ref{vac} number $p$ does not divide $|R_i|$. Thus $G/G_{i-1}$ possesses an element $g\not\in R_i$ of order $p$, which acts by conjugation on the group $R_i$. According to Lemma \ref{vac2} we have $R_i<C_{G/G_{i-1}}(g)$, but $\pi(C_{G/G_{i-1}}(g))$ has no numbers from $\Pi$. Hence for every $1\leq i< m$ the intersection $\pi(R_i)\cap \Pi$ contains at most one number. Suppose that there exist two numbers $1\leq h<l<m$ such that $|R_h|$ and $|R_l|$ are divided by distinct numbers from $\Pi$. The group $G$ contains a normal series $1 <G_h < G_l<G_m$ and numbers $r,\ q,\ p\in \Pi$ such that $r$ divides $|G_h|$, $q$ divides $|G_l/G_h|$ and $p$ divides $|G_m/G_l|$. The numbers $r,\ q,\ p$ are pairwise non-adjacent in $GK(G)$ which contradicts Lemma \ref{vas2}. The lemma is proved.
\smallskip

Let us complete the proof of Theorem 1. Lemmas \ref{tc} and \ref{pi} imply that for every $n\geq26$, $n\neq 27$, there exists $r\in \overline{\Pi}$ such that $n-r\not\in\{3,4,6,10\}$. We show that under $n=27$ this statement is also true. Since $\Pi=\{17,19\}$, we have to prove that $19\in\overline \Pi$. Suppose that $19$ divides $|G_{m-1}|$. We have $R\simeq A_k$, where $k\in[23, 27]$. The group $R$ contains the Frobenius group $F$ with the kernel of order $23$ and the compliment of order $11$. By Lemma \ref{vas2} the full preimage of the group $F$ contains an element of order $11\cdot 19$ or $23\cdot 19$; contradiction.

Since a Sylow $r$-subgroup of the group $R$ has order $r$ and $r$ does not divide $|G_{m-1}|\cdot|G/G_m|$, then a Sylow $r$-subgroup of the group $G$ has order $r$ also. Suppose that $g\in G$ is an element of order $r$, $\overline{g}$ is the image of $g$ in $R$. According to Lemma \ref{centralizer2} we have $C_G(g)=\langle g\rangle\times M$, $\mu(M)=\mu(A_{n-r})$. The group $A_s$ is recognizable for every $5\leq s<n$, $s\neq 6, 10$ which implies that $M\simeq A_{n-r}$. Suppose that $C_i=C_G(g)\cap
G_i$, $0\leq i \leq t$. Since $n-r\geq5$, the group $M$ is simple, hence there exists $i$ such that $M\leq C_i$ and
$M\cap C_{i-1}=1$. By Lemma \ref{factor} we have $C_R(\overline{g})\simeq C_{G}(g)/(G_{m-1}\cap C_{G}(g))$. Since
$C_R(\overline{g})\neq \langle \overline{g}\rangle$, then $i=m$. By Lemma \ref{sp} we have $k=n$. Lemmas \ref{ASS} and \ref{zav2} imply that $G\simeq L$. Theorem 1 is proved.

~

Gorshkov Ilya Borisovich

4 Acad. Koptyug avenue, 630090 Novosibirsk Russia

Sobolev Institute of Mathematics,

ilygor@ngs.ru,

(383)3634613

\bigskip
\begin{thebibliography}{1}
\bibitem{MazKon} A. S. Kondrat'ev, V. D. Mazurov, "Recognition of alternating groups of prime degree from their element orders", Siberian Mathematical Journal, 41:2 (2000), 359–369.
\bibitem{Zav}  A. V. Zavarnitsin, "Recognition of Alternating Groups of Degrees r+1 and r+2 for Prime r and the Group of Degree 16 by Their Element Order Sets", Algebra and Logic, 39:6 (2000), 648–661.
\bibitem{Shi} R. Brandl, W. Shi, "Finite groups whose element orders are consecutive integers", J.Algebra, 143:2 (1991), 388-400.
\bibitem{Maz2} V. D. Mazurov, "Characterizations of finite groups by sets of element orders", Algebra and Logic, 37:6 (1998), 651–666.
\bibitem{Changg} Changguo Shao, Qinhui Jiang, "A new characterization of $A_{22}$ by its spectrum", Communications in algebra, 38:6 (2010), 2138-2141.
\bibitem{ZavMaz} A. V. Zavarnitsin, V. D. Mazurov, " Element orders in coverings of symmetric and alternating groups", Algebra and Logic, 38:3 (1999), 296-315.
\bibitem{Vac} I. A. Vakula, "On the structure of finite groups isospectral to an alternating group", Proceedings of the Steklov Institute of Mathematics, 16:3 (2010), 45–60.
\bibitem{VasVd} A. V. Vasiliev, E. P. Vdovin, "An Adjacency Criterion for the Prime Graph of a Finite Simple Group", Algebra and Logic, 44:6 (2005), 682–725.
\bibitem{Khu} E. I. Khukhro, Nilpotent groups and their automorphisms, De Gruyter, Berlin, 1993.
\bibitem{Vas} A. V. Vasil'ev, "On connection between the structure of a finite group and the properties of Its prime graph", Siberian Mathematical Journal, 46:3 (2005), 511–522.

\end{thebibliography}
\end{document}